
\baselineskip=14pt
\parskip=10pt

\magnification=\magstephalf

\def\W{{\cal W}}
\def\1{{\overline{1}}}
\def\2{{\overline{2}}}
\parindent=0pt
\overfullrule=0in

\def\frac#1#2{{#1 \over #2}}
\centerline
{\bf On Invariance Properties of Entries of Matrix Powers}
\bigskip
\centerline
{\it Shalosh B. EKHAD and Doron ZEILBERGER}

\bigskip

{\bf Abstract}: A  few years ago, Peter Larcombe discovered an amazing property  regarding two by two matrices. For any such $2 \times 2$ matrix $A$,
the ratios of the two anti-diagonal entries is the same for all powers of $A$.
We discuss extensions to higher dimensions, and give a short bijective proof of Larcombe and Eric Fennessey's elegant extension
to tri-diagonal matrices of {\it arbitrary} dimension. This article is accompanied by a Maple package.

{\bf Peter Larcombe's Surprising Discovery}

In [2], Peter Larcombe gave four proofs of a seemingly new amazing property of a $2 \times 2$ matrix, for any such matrix
(we denote the $(i,j)$ entry of a matrix $B$ by $B_{ij}$)

$$
A_{12} \cdot (A^m)_{21}=A_{21} \cdot (A^m)_{12} \quad,
$$
for {\it all} positive integers $m$.

We first observe that, in hindsight (but only in hindsight!) this is not that surprising. More generally, for a general $n \times n$ matrix $A$,
and any subset $S$ of cardinality $n+1$ 
of the set of $n^2$ entries $\{(i,j) \,|\, 1\leq i,j \leq n\}$, there exist polynomials $q_s(A)$ in the entries of $A$ (independent of $m$) such that
$$
\sum_{s \in S} \, q_s \cdot (A^m)_s \, = \, 0 \quad,
\eqno(1)
$$
for all $m>1$.

This fact follows from the {\bf Cayley-Hamilton} equation that says that
If $P_A(x):=\det(A-x\,I)$ is the {\bf characteristic polynomial} of $A$, then the $n \times n$ matrix $P_A(A)$ equals the {\bf zero matrix} ${\bf 0}$ . Writing
$$
P_A(X)=\sum_{k=0}^{n} p_k x^k \quad ,
$$
we have
$$
\sum_{k=0}^{n} p_k A^k= {\bf 0} \quad,
$$
where ${\bf 0}$ is the all-zero matrix.
Multiplying by $A^m$ we get
$$
\sum_{k=0}^{n} p_k A^{m+k} = {\bf 0} \quad,
$$
for all $m$.

Taking the $ij$ entry, we have that each of the $n^2$ sequences $(A^m)_{ij}$ satisfy the {\bf same} $n^{th}$-order linear recurrence equation with {\bf constant coefficients}

$$
\sum_{k=0}^{n} p_k (A^{m+k})_{ij} = 0 \quad .
$$

It is well-known and easy to see ([1][8]) that any $n+1$ sequences that satisfy the {\bf same} recurrence of order $n$, must be {\bf linearly dependent}.
Also, in order to find the relation, it is enough to find  $n+1$ values of the $q_s$ for which $(1)$ holds for the  $m \,= \, 1, \, 2, \dots, n$. Then this linear combination
also satisfies that very same linear recurrence, and since it {\bf vanishes} at the first $n$ {\bf initial conditions} it must be {\bf identically zero}.

If our subset $S$ of entries only consists of non-diagonal entries, we can do even better, Eq. $(1)$ is true for any set $S$ of $n$ non-diagonal entries.

Note that the Cayley-Hamilton equation implies that

$$
\sum_{k=1}^{n} p_k A^k \quad
$$
is a {\bf diagonal} matrix (namely $-\det(A) \bf I$), hence for a non-diagonal entry $ij$ ($i \neq j$), $(A^m)_{ij}$ always satisfies the {\bf same} linear recurrence equation (with constant coefficients)
of order $n-1$, hence any $n$ such non-diagonal entries must be {\bf linearly dependent}. In the original case of a $2 \times 2$ matrix discussed in [2], it follows that the $(1,2)$ and $(2,1)$ entries of $A^m$ 
always satisfy the {\bf same} relation as those of $A$, hence we have yet-another-proof (without equations!) of Larcombe's amazing discovery.

But what about higher dimensions? Now things get much more complicated, and we need a computer algebra system (in our case Maple). For example, we have the
following

{\bf Theorem}: Let $A=(a_{ij})_{1\leq i,j\leq 3}$ be a $3 \times 3$ matrix, then for {\bf all} $m \geq 1$, we have
$$
 \left(a_{1,2} a_{2,1} a_{2,3}-a_{1,3} a_{2,1} a_{2,2}+a_{1,3} a_{2,1} a_{3,3}-a_{1,3} a_{2,3} a_{3,1}\right) \cdot (A^m)_{12}
$$
$$
+\left(a_{1,2} a_{2,3} a_{3,1}-a_{1,3} a_{2,1} a_{3,2}\right) \cdot (A^m)_{13}
$$
$$
+\left(-a_{1,2}^{2} a_{2,3}+a_{1,2} a_{1,3} a_{2,2}-a_{1,2} a_{1,3} a_{3,3}+a_{1,3}^{2} a_{3,2}\right) \cdot (A^m)_{21} \,=\, 0 \quad .
$$

There are three more such theorems (up to trivial isomorphism)  for the $n=3$ case, while 
there are $27$ inequivalent cases for $n=4$. They can all be found in the following output file

{\tt https://sites.math.rutgers.edu/\~{}zeilberg/tokhniot/oLarcombe1.txt} \quad .

We were unable to find the corresponding relations for $n=5$, they got too complicated!

\vfill\eject

{\bf A bijective proof of the Larcombe-Fennessey theorem about Tridiagonal matrices}

Any matrix identity for a fixed dimension is essentially {\it high-school algebra} and can be verified by a computer algebra system, even if the power $m$ is arbitrary. But the following theorem
of Larcombe and Fennessey, regarding {\bf tridiagonal} matrices of {\it arbitrary} dimension is {\it university algebra} and is more interesting.

{\bf Theorem} (Larcombe and Fennessey [6]) : Let $A$ be a general tridiagonal $n \times n$ (for {\it any} $n \geq 2$), then for all $1 \leq i<n$, and all $m \geq 1$, we have
$$
a_{i,i+1} \cdot (A^m)_{i+1,i} \, = \,  a_{i+1,i} \cdot (A^m)_{i,i+1} \quad .
\eqno(2)
$$

We will give a {\bf combinatorial proof}. Fix $n$ and let $a_{i,j}$ ($1 \leq i,j \leq n$) be $n^2$ {\bf commuting} indeterminates. It follows immediately from the definition of matrix multiplication that
the $(i,j)$ entry of $A^m$ is the {\bf weight-enumerator} of the set of $(m+1)$-letter words in the alphabet $\{1,2,\dots,n\}$ whose first letter is $i$  and  last letter is $j$, with the weight
$$
Weight(w_1 \dots w_{m+1}):= a_{w_1\,w_2} \, \cdot \, a_{w_2\,w_3}\,  \cdots \, a_{w_{m-1}\,w_m} \, \cdot\,  a_{w_{m}\,w_{m+1}}  \quad .
$$

If our matrix is {\it tridiagonal}, then all the words are {\it continuous} i.e. after the letter $i$ can only come one of the (up to) three letters $\{ i-1,i,i+1 \}$. For example
if $n=5$ then the following is a legal word

$$
2333234333221122112234455443 \quad,
$$
but the following one is {\bf not}
$$
233312 \quad,
$$
because after the fourth letter, that is a `$3$', comes the letter `$1$'.

Let $\W_{m}(i,j)$ be the set of legal $(m+1)$-letter words in the alphabet $\{1,2, \dots, n\}$  that start with the letter $i$ and end with the letter $j$.  Its weight-enumerator (i.e. sum of the weights of its members) is
$(A^{m})_{ij} $, where  now $A$ is a generic $n \times n$ {\bf tridiagonal} matrix. 
For ease of type-setting let $i':=i+1$.

Note that:

$\bullet$ The left side of $(2)$ is is the weight-enumerator of the set of words, $i \, \W_{m}(i',i) $, which is the set of {\bf legal} $(m+2)$-letter words that start and end with the letter $i$, and whose {\bf second} letter if $i'$.

$\bullet$ The right side of $(2)$ is is the weight-enumerator of $i'\, \W_{m}(i,i') $ which is the set of {\bf legal} $(m+2)$-letter words that start and end with the letter $i'$, and whose {\bf second} letter if $i$.

We claim that The mapping 

$$
T_i\,:\, i \,\W_{m}(i',i) \rightarrow  i' \,\W_{m}(i,i') \quad,
$$ 
to be defined next, is a {\bf weight-preserving} bijection.

Let $w=w_1 w_2 \dots w_{m+2}$ be a member of  $i \, \W_{m}(i',i)$, then of course $w_1=i$ and $w_2=i'$, and $w_{m+2}=i$. Let $k$ be the {\bf smallest} index such that $w_k=i', w_{k+1}=i$. Of course it exists (by ``{\it continuity}'').

{\bf Case I}: $k=2$. 

If $m=1$ then the word must be $i\,i\,'i$ and we map it to $i' \, i \, i'$. 

Otherwise we can write 

$$
w\,=\, i\,i'\,i\,u\,i  \quad,
$$ 

for some  $(m-2)$-letter word $u$, and we define 

$$
T_i(w):=\, i' \, i \, u \, i \, i' \quad,
$$

that of course belongs to  $i' \, \W_{m}(i,i') $.

{\bf Case II}: $k=m+1$, then we can write  
$$
w\,=\, i \, i' \, u \, i' \, i \quad ,
$$ 
for some $(m-2)$-letter word $u$, and we map it to the

$$
T_i(w):= \,i' \, i \, i' \, u \,i'  \quad,
$$
that of course belongs to  $i' \, \W_{m}(i,i') $.

{\bf Case III}: $2<k<m+1$. Then we can write 

$$
w\, =\,  i\,i' \, u \, i' \, i \, v \, i \quad,
$$ 
for some words $u$ and $v$, whose total length is $m-2$, and we define
$$
T_i(w):=\, i' \, i \, v \, i \, i' \, u \, i'  \quad,
$$

that of course belongs to  $i'\, \W_{m}(i,i') $.

Let's state the inverse mapping 
$$
U_i \,:\,  i' \, \W_{m}(i,i') \rightarrow  i \, \W_{m}(i',i) \quad.
$$ 

Let $w=w_1 w_2 \dots w_{m+2}$ be a member of  $i' \, \W_{m}(i,i')$, then of course $w_1=i'$ and $w_2=i$, and $w_{m+2}=i'$. Let $k$ be the {\bf largest} index such that $w_k=i, w_{k+1}=i'$. Of course it exists (by ``{\it continuity}'').

{\bf Case I}: $k=m+1$. If $m=2$ the $w=\,i' \, i \, i'$ and we let $U_i(w)$ be $i\, i' \, i$. Otherwise we can write 

$$
w\,=\, i' \, i \, u \, i \, i' \quad,
$$ 

for some $(m-2)$-letter word $u$, and we define 

$$
U_i(w):= \, i \,i' \, i \, u \, i \quad .
$$

{\bf Case II}: $k=2$. We can write 

$$
w\,=\, i' \, i \, i' \, u \, i' \quad
$$ 

and we define 

$$
U_i(w):= \, i \, i' \, u \, i' \, i  \quad .
$$

{\bf Case III}: $2<k<m+1$. We can write 

$$
w\,=\, i' \, i \, v \, i \, i' \, u \, i' \quad ,
$$ 

for some words $u$ and $v$ whose lengths add-up to $m-2$, and we define 

$$
U_i(w):=\, i \, i' \, u \, i' \, i \, v \, i \quad .
$$

Readers are welcome to play with the Maple package

{\tt https://sites.math.rutgers.edu/\~{}zeilberg/tokhniot/Larcombe.txt} \quad,

that contains procedure $\tt REL$ to discover generalized Larcombe relations (mentioned above) and also implements the above bijection
(procedures {\tt Ti} and {\tt Ui}, and {\tt CheckTi} verifies it empirically).

{\bf References}

[1] Manuel Kauers  and Peter Paule, {\it ``The Concrete Tetrahedron''}, Springer, 2011.

[2] P. J. Larcombe, {\it A note on the invariance of the general $2 \times 2$  matrix anti-diagonal ratio with increasing matrix powers: four proofs},
Fib. Quarterly {\bf 53}, 360-364 (2015).

[3] P. J. Larcombe and E. J. Fennessey, {\it A new tri-diagonal matrix invariance property}, Palest. J. Math. {\bf 7}, 9-13 (2018).

[4] P. J. Larcombe and E. J. Fennessey, {\it A note on two rational invariants for a particular $2 \times 2$ matrix}, Palest. J. Math. {\bf 7}, 410–413 (2018).

[5] P. J. Larcombe and E. J. Fennessey, {\it On anti-diagonals ratio invariance with exponentiation of a $2 \times 2$ matrix: two new proofs}, Palest. J. Math. {\bf 9}, 12–14 (2020).

[6] P. J. Larcombe and E. J. Fennessey, 
{\it A formalised inductive approach to establish the invariance of antidiagonal ratios with exponentiation for a tri-diagonal matrix of fixed dimension}, 
Palest. J. Math. {\bf 9}, 670-672 (2020).

[7] P. J. Larcombe and E. J. Fennessey, 
{\it A short graph-theoretic proof  of the  $2 \times 2$ matrix anti-diagonal ratio invariance with exponentiation}, Palest. J. Math. 10, 102-103 (2021).

[8] Doron Zeilberger, {\it The C-finite Ansatz}, Ramanujan Journal {\bf 31} (2013), 23-32.\hfill\break
{\tt https://sites.math.rutgers.edu/\~{}zeilberg/mamarim/mamarimhtml/cfinite.html} \quad .

\bigskip
\hrule
\bigskip
Shalosh B. Ekhad and Doron Zeilberger, Department of Mathematics, Rutgers University (New Brunswick), Hill Center-Busch Campus, 110 Frelinghuysen
Rd., Piscataway, NJ 08854-8019, USA. \hfill\break
Email: {\tt [ShaloshBEkhad, DoronZeil] at gmail dot com}   \quad .

Written: {\bf July 27, 2021}.

\end